\documentclass{svproc}

\usepackage{amssymb}
\usepackage{amsfonts}
\usepackage{mathtools}
\usepackage{mathrsfs}
\usepackage{bm}
\usepackage{bbm}
\usepackage{caption}
\usepackage{algorithm,algorithmicx,algpseudocode}
\usepackage{graphicx}
\usepackage{todonotes}
\usepackage{setspace}
\usepackage{mathtools}
\usepackage{url}

\newcommand{\bE}{{\mathbf E}}
\newcommand{\E}{{\mathbf E}}

\newcommand{\wrt}{with respect to }

\newcommand{\R}{{\mathbb R}}

\newcommand{\cE}{{\mathcal E}}

\begin{document}
\mainmatter              
\title{Adaptive importance sampling with forward-backward stochastic differential equations}
\titlerunning{Importance sampling and FBSDE}  
%
\author{Omar Kebiri\inst{1,2} \and Lara Neureither\inst{2} \and Carsten Hartmann\inst{2}}
%
%
 \tocauthor{Lara Neureither and Carsten Hartmann}

\institute{$^{1}$Laboratory of Statistics and Random Modeling\\ University of Abou Bekr Belkaid, Tlemcen, Algeria\\[1ex]
$^{2}$Brandenburgische Technische Universit\"at Cottbus-Senftenberg\\ Cottbus, Germany}
\maketitle              

\begin{abstract}
	 We describe an adaptive importance sampling algorithm for rare events that is based on a dual stochastic control formulation of a path sampling problem. Specifically, we focus on path functionals that have the form of cumulate generating functions, which appear relevant in the context of, e.g.~molecular dynamics, and we discuss the construction of an optimal (i.e. minimum variance) change of measure by solving a stochastic control problem. We show that the associated semi-linear dynamic programming equations admit an equivalent formulation as a system of uncoupled forward-backward stochastic differential equations that can be solved efficiently by a least squares Monte Carlo algorithm. We illustrate the approach with a suitable numerical example and discuss the extension of the algorithm to high-dimensional systems.

\end{abstract}

\section{Introduction}

The simulation of rare events is among the key challenges in computational statistical mechanics which involves fields such as molecular dynamics \cite{HS12}, material science \cite{ER93} or climate modelling \cite{WB16}. Concrete examples include the study of critical phase transitions in many-particle systems or the estimation of small transition probabilities in protein folding. Estimating small probabilities by Monte Carlo is tricky, because the standard deviation of the corresponding statistical estimator is typically larger than the quantity to be estimated.
One technique to improve the efficiency of estimators for
small probabilities is importance sampling. Here the idea is to sample from another distribution under which the rare event is no longer rare and then correct (i.e.~reweight) the estimator with the appropriate likelihood ratio. Designing such a change of measure so that the variance of the reweighted estimator stays bounded is not at all a trivial task, and several methods have been developed to cope with this issue; for an overview, we refer to the standard textbooks \cite{AssGlynn,Rubin08} and the references therein.

Here we consider adaptive importance sampling strategies where the change of measure is mediated by an exponential tilting of the reference probability measure. For stochastic differential equations, this exponential tilting can be interpreted as a control that changes the drift of the stochastic dynamics. Adaptive importance sampling has been predominantly studied in the context of small noise diffusions, for which the optimal control can be computed from the zero viscosity limit of the corresponding dynamic programming equation \cite{DW04,DW07,VW12}. In this case the value function of the zero viscosity (deterministic) control problem is equal to the large deviations rate function that describes the exponential tails of the rare events under consideration, and as a consequence, the change of measure captures the rare events statistics and results in estimators that, under certain assumptions,  have uniformly bounded relative error.

Here we follow a different route, in that we do not resort to large deviations asymptotics, but rather try to compute the zero-variance change of measure from a suitable approximation of the dynamic programming equation that is underlying the stochastic control problem; in contrast to our previous works \cite{HS12,ZWHWS14}, in which the change of measure has been obtained by solving the corresponding variational problem directly, we here focus on the reformulation of the underlying dynamic programming equation as a system of forward-backward stochastic differential equations (see, e.g.~\cite{peng93,to}) that is solved by a least squares Monte Carlo algorithm \cite{BenderSteiner2012,GobetEtal2005}. Our approach is partly inspired by related duality techniques in financial mathematics \cite{HK04,Rogers02}, but exploits the specific duality structure of the change of measure problem; see also \cite{P09} for a survey of related approaches in financial mathematics.

The paper is organised as follows: In Section \ref{sec:is} we introduce our stochastic dynamics, the corresponding path space free energy and its dual variational characterisation.    
Section \ref{sec:fbsde} deals with the formulation of the free energy sampling problem and the (dual) optimal control problem as a forward-backward stochastic differential equation (FBSDE, in short). The numerical solution of the FBSDE that can be used to either directly compute the free energy or to approximate the optimal control that generates the minimum variance importance sampling scheme is the topic of Section \ref{sec:mc}, with a simple numerical illustration presented in Section \ref{sec:num}. The article concludes in Section \ref{sec:out} with a short summary and a discussion of open problems and future work.

\section{Importance sampling in path space}\label{sec:is}

Let $X=(X_s)_{s\geqslant 0}$ be the solution of
\begin{equation}\label{sde}
dX_{s} = b(X_{s},s)\, ds + \sigma(X_{s}) dB_{s} \,,\quad X_0=x\,,
\end{equation}
where $X_s \in \mathbb{R}^d$, $b$ and $\sigma$ are smooth
 drift and noise coefficients, and $B$ is an $m$-dimensional standard Brownian motion where in general $m\leqslant d$. Our standard example will be a non-degenerate diffusion in an energy landscape,
 \begin{equation}\label{sdegrad}
dX_{s} = -\nabla U(X_{s})ds + \sigma dB_{s} \,,\quad X_0=x\,,
\end{equation}
 with smooth potential energy function $U$ and $\sigma>0$ constant.
We assume throughout this paper that the functions $b,\sigma, U$ are such that Equations (\ref{sde}) or \eqref{sdegrad} have unique strong solutions for all $s\geqslant 0$.
Now let $W$ be a continuous functional
\begin{equation}\label{Wfunc}
W_\tau(X) = \int_{0}^{\tau} f(X_s,s)\,ds + g(X_\tau)\,,
\end{equation}
of $X$ up to some bounded stopping time $\tau$ where $f,g$ are bounded and sufficiently smooth, real valued functions. 

 \begin{definition}[Path space free energy] Let $X$ be the solution
  of Equation (\ref{sde}) and let $W_\tau=W_{\tau}(X)$ be defined by Equation (\ref{Wfunc}). The quantity
	\begin{equation}\label{cgf}
	 \gamma = -\log \bE\left[\exp(-W_\tau)\right]
	\end{equation}
	is called the free energy of $W_\tau$ where the expectation is understood \wrt the realisations of (\ref{sde}) for given a initial condition $X_{0}=x$.
\end{definition}

\subsection{Donsker--Varadhan variational formula for the free energy}

The adaptive importance sampling strategy described below is based on a variational characterization of (\ref{cgf}) in terms of a change of measure.  To make it precise, we define $P$ to be the probability measure on the space $\Omega = C([0,\infty),\R^n)$ of continuous trajectories that is induced by the Brownian motion $B$ in (\ref{sde}). We denote the expectation \wrt $P$ by $\bE[\cdot]$.
In abstract form, the Donsker--Varadhan variational principle \cite{Ellis1985} states
\begin{equation}\label{dv}
\gamma = \inf_{Q\ll P} \left\{\bE_Q[W_{\tau}] + D(Q|P)\right\}\,,
\end{equation}
where $Q\ll P$ stands for {absolute continuity} of $Q$ \wrt $P$, and
\begin{equation}\label{kl}
D(Q|P) = \begin{cases}
\displaystyle\int_\Omega \log\frac{dQ}{dP}(\omega) \,dQ(\omega) & \textrm{if } Q\ll P\\
+\infty & \textrm{else}\,.
\end{cases}
\end{equation}
denotes the \emph{relative entropy} or \emph{Kullback--Leibler divergence} between $Q$ and $P$. Note that $D(Q|P)=\infty$ when $Q$ is not absolutely continuous \wrt $P$, therefore it is sufficient to take the infimum in (\ref{dv}) over all path space measures $Q\ll P$. If $W_\tau\geqslant 0$, it is a simple convexity argument (see, e.g.,~\cite{DMR96}), which shows that the minimum in Equation (\ref{dv}) is attained at $Q^*$ given by
\begin{equation}\label{mindv}
\frac{dQ^*}{dP}\bigg|_{\mathcal{F}_{\tau}} = \exp(\gamma - W_\tau)\,,
\end{equation}
where $\varphi|_{\mathcal{F}_{\tau}}$ denotes the restriction of the path space density $\varphi=dQ^*/dP$ to the $\sigma$-algebra ${\mathcal{F}_{\tau}}\subset\cE$ that is generated by the Brownian motion $B$ up to time $\tau$.\footnote{More precisely, $\varphi|_{\mathcal{F}_{\tau}}$ is understood as the restriction of the measure $Q^*$ defined by $dQ^*=\varphi dP$ to the $\sigma$-algebra ${\mathcal{F}_{\tau}}$ that contains all measurable sets $E\in\cE$, with the property that for every $t\geqslant 0$ the set $E\cap\{\tau\leqslant t\}$ is an element of the $\sigma$-algebra $\mathcal{F}_{t}=\sigma(X_s:0\leqslant s\leqslant t)$ that is generated by all trajectories $(X_s)_{0\leqslant s\leqslant t}$ of length $t$.}
By the strict convexity of the exponential function, it holds that $Q^{*}$-a.s. \cite{HartmannEnt2017}
\begin{equation}\label{zerovar1}
\bE\left[\exp(-W_\tau)\right] = \exp(-W_\tau)\left(\frac{dQ^*}{dP}\bigg|_{\mathcal{F}_{\tau}}\right)^{-1}
\end{equation}
or, equivalently,
\begin{equation}\label{zerovar2}
\gamma = W_\tau + \log\left(\frac{dQ^*}{dP}\bigg|_{\mathcal{F}_{\tau}}\right)\,.
\end{equation}
That is, $Q^{*}$ defines a zero-variance change of measure. (Note that the inverse of the Radon--Nikodym derivative in (\ref{zerovar1}) exists since $W_\tau$ is bounded.)

\subsection{Related stochastic control problem}

The only admissible change of measure from $P$ to $Q$ such that $D(Q|P)<\infty$ results in a change of the drift in Equation (\ref{sde}). Specifically, let $u$ be a process with values in $\R^{m}$ that is adapted to $B$ and that satisfies
\begin{equation}\label{novikov}
\bE\left[\exp\left(\frac{1}{2}\int_{0}^{\tau}|u_{s}|^{2}\,ds\right)\right] < \infty\,.
\end{equation}
Further define the auxiliary process
\begin{equation*}
B^u_{t} = B_{t} - \int_{0}^{t}u_{s}\,ds\,,
\end{equation*}
so that (\ref{sde}) can be expressed as
\begin{equation}\label{drivenSDE}
dX_{s} = \left(b(X_{s},s) + \sigma(X_{s})u_{s}\right)ds + \sigma(X_{s}) dB^u_{s} \,,\quad X_{0}=x\,.
\end{equation}
By construction, $B^u$ is not a Brownian motion under
$P$, but by Girsanov's Theorem (see, e.g.,~\cite{O10}, Theorem~8.6.4) there exists a measure $Q$ defined by
\begin{equation}\label{girsanov}
\frac{dQ}{dP}\bigg|_{\mathcal{F}_{\tau}}= \exp\left(\int_{0}^{t}u_{s}\cdot dB^u_{s} + \frac{1}{2} \int_{0}^{t}|u_{s}|^{2}\,ds\right) \,
\end{equation}
so that $B^u$ is a standard Brownian motion under $Q$. (The Novikov condition (\ref{novikov}) guarantees that $Q$ is a probability measure.) Inserting (\ref{girsanov}) into  (\ref{dv}), using that $B^u$ is a Brownian motion \wrt $Q$, it follows that (cf.~\cite{BD98,DMR96}):
\begin{equation}
\gamma = \inf_{u}\bE_Q\!\left[ \int_{0}^{\tau}f(X_{s},s) + \frac{1}{2} |u_{s}|^{2}\,ds + g(X_\tau)\right]\,,
\end{equation}
with $X$ being the solution of Equation (\ref{drivenSDE}). Since the distribution of $B^u$ under $Q$ is the same as the distribution of $B$ under $P$, an equivalent representation of the last equation is
\begin{equation}\label{dv2}
\gamma = \inf_{u}\bE\!\left[ \int_{0}^{\tau}f(X^u_{s},s) + \frac{1}{2} |u_{s}|^{2}\,ds + g(X^u_\tau)\right]\,.
\end{equation}
where $X^u$ is the solution of the controlled SDE
\begin{equation}\label{drivenSDE2}
dX^u_{s} = \left(b(X^u_{s},s)+\sigma(X^u_s)u_s\right)ds + \sigma(X^u_{s}) dB_{s} \,,\quad X^u_0=x\,,
\end{equation}
with $B$ being a standard $m$-dimensional Brownian motion under the probability measure $P$. The Donsker--Varadhan variational principle (\ref{dv}) and zero-variance property (\ref{zerovar1}) of the probability measure $Q^{*}$, for which equality in (\ref{dv}) is attained, have the following stochastic control analogue (see \cite[Thm.~3.1]{HartmannEnt2017}):
\begin{theorem}\label{thm:zerovar}
	Let $T>0$ and $\tau_O = \inf\{s>0\colon X^u_{s}\notin O\}$ for an open and bounded set  $O\subset\R^{n}$ with smooth boundary $\partial O$. Further define $\tau=\tau_O\wedge T$ and    	
	\begin{equation}\label{mgf}
	 \Psi(x,t)
	 =\bE\!\left[\exp\left(-\int_{t}^{\tau}
	 f(X_s,s)\,ds - g(X_\tau)\right)\Bigg| X_t=x\right]
	\end{equation}
	as the exponential of the negative free energy, considered as a function of the initial condition $X_t=x$ with $0\leqslant t\leqslant \tau\leqslant T$. Then, the path space measure $Q^*$ induced by the feedback control
	\begin{equation}\label{oc}
	u^*_s = \sigma(X^{u^*}_s)^T\nabla_x \log\Psi(X^{u^*}_s,s)
	\end{equation}
	and (\ref{drivenSDE2}) yields a zero variance estimator, i.e.,
	\begin{equation}\label{zeroSDE}
		\Psi(x,0) = \exp\left(-\int_{0}^{\tau}
	  f(X^{u^{*}}_s,s)\,ds - g(X^{{u^{*}}}_\tau)\right)\left(\frac{dQ^*}{dP}\bigg|_{\mathcal{F}_{\tau}}\right)^{-1}\;\textrm{$Q^*$-a.s.}
	\end{equation}
\end{theorem}

\section{From dynamic programming to forward-backward SDE}\label{sec:fbsde}

Following the route taken by \cite{BD98}, it can be shown that the control $u^*$ in (\ref{oc}) is the unique minimiser of the following stochastic control problem: minimise
\begin{equation}\label{cost}
J(u) = \bE\!\left[\int_{0}^{\tau}
f(X^u_s,s) + \frac{1}{2}|u_s|^2\,ds + g(X^u_\tau)\right]
\end{equation}
over all measurable and square integrable Markovian controls $u$, such that the controlled SDE (\ref{drivenSDE2}) has a unique strong solution. Now let
\begin{equation}\label{cost2go}
V(x,t) = \min_u\bE\!\left[\int_{t}^{\tau}
f(X^u_s,s) + \frac{1}{2}|u_s|^2\,ds + g(X^u_\tau)\bigg|\, X_t^u=x\right]
\end{equation}
be the associated value function (or: optimal cost-to-go). Further define $E=O\times [0,T)$ and let $\partial E^{+} = \left(\partial O\times [0,T)\right)\cup\left(O\times \{T\}\right)$ be the terminal set of the augmented process $(X^{u}_{s},s)_{s\geqslant 0}$, such that $\tau=\tau_{O}\wedge T$ can be recast as
\begin{equation}
\tau = \inf\!\left\{s>0\colon (X^{u}_{s},s)\notin E\right\}\,.
\end{equation}
Assuming sufficient regularity of the coefficients $b,\sigma, f,g$ and $\partial O$, a necessary and sufficient condition for $u=u^{*}$ being optimal is that   (see \cite[Sec. VI.5]{FM06})
\begin{equation}\label{oc2}
	u^*_s = - \sigma(X^{u^*}_s)^T\nabla_x V(X^{u^*}_s,s)
\end{equation}
where $V\in C^{2,1}(E)\cap C(\partial E^{+})$ solves the dynamic programming equation 
\begin{equation}\label{hjb}
\begin{aligned}
\partial_{t} V + LV + h(s,x,V,\sigma^{T}\nabla_{x}V) & = 0\quad \textrm{in $E$}\\ V & = g\quad \textrm{on $\partial E^{+}$}\,,
\end{aligned}
\end{equation}
with nonlinearity
\begin{equation}\label{F}
h(s,x,y,z) = - \frac{1}{2} |z|^{2} + f(x,s)
\end{equation}
and the infinitesimal generator of the control-free process $X_{t}$,
\begin{equation}\label{L}
L = \frac{1}{2}\sigma\sigma^{T}\colon\nabla^{2}_{x} + b\cdot \nabla_{x}\,.
\end{equation}
For the derivation of (\ref{oc2})--(\ref{hjb}) from the Feynman--Kac representation formula for the free energy (\ref{cgf}), we refer to  \cite[Sec.~6]{HartmannEntropy2014}.

\subsection{FBSDE representation of the dynamic programming equation}

We will now recast the semi-linear, parabolic boundary value problem for $V\in C^{2,1}(E)\cap C(\partial E^{+})$. To this end, define the processes
\begin{equation}\label{YandZ}
Y_s = V(X_s,s)\,,\quad Z_s = \sigma(X_s)^T\nabla_x V(X_s,s)\,
\end{equation}
with $X$ denoting the solution of the uncontrolled SDE (\ref{sde}) with infinitesimal generator (\ref{L}). Applying Ito's formula to $Y$, using that $V$ is a classical solution to (\ref{hjb}), we obtain the following backward SDE (BSDE)
\begin{equation}\label{bsde}
dY_s = -h(s,X_s,Y_s,Z_s)ds + Z_s\cdot dB_s\,,\quad Y_\tau = g(X_\tau)\,
\end{equation}
for the pair $(Y,Z)$. Note that, by definition, $Y$ is continuous and adapted to $X$, and $Z$ is predictable and a.s.~square integrable, i.e.,
\begin{equation}
\int_0^\tau |Z_s|^2\,ds <\infty\,,
\end{equation}
in accordance with the interpretation of $Z_s$ as a control variable. Further note that (\ref{bsde}) must be understood as a \emph{backward} SDE rather than a \emph{time-reversed} SDE, since, by definition, $Y_s$ at time $s<\tau$ is measurable with respect to the filtration generated by the Brownian motion $(B_r)_{0\leqslant r\leqslant s}$, whereas a time-reversed version of $Y_s$ would depend on $B_\tau$ via the terminal condition $Y_\tau=g(X_\tau)$, which would require a larger filtration.

By exploiting the specific form of the nonlinearity (\ref{F}) that appears as the driver $h$ in the backward SDE (\ref{bsde}) and the fact that the forward process $X$ is independent of $(Y,Z)$, we obtain the following representation of the solution to the dynamic programming equation (\ref{hjb}):
\begin{equation}\label{fbsde}
\begin{aligned}
dX_s & = b(X_s,s)ds + \sigma(X_s)\,dB_s\,,\quad X_t=x\\
dY_s & = -f(X_s,s)ds + \frac{1}{2}|Z_s|^2 + Z_s\cdot dB_s\,,\quad Y_\tau = g(X_\tau)\,.
\end{aligned}
\end{equation}
The solution to (\ref{fbsde}) now is a triplet $(X,Y,Z)$, and since $Y$ is adapted, it follows that $Y_t$ is a deterministic function of the initial data $(x,t)$ only. Since $g$ is bounded, the results in \cite{Kobylanski2000} entail existence and uniqueness of (\ref{bsde}); see also \cite{bgm,BKM}. 
As a consequence (see e.g.~\cite{PP92} or \cite[Prop.~3.1]{BBM86}),
\begin{equation}\label{bsdevalue}
Y_t = V(x,t)\quad \textrm{(a.s.)}
\end{equation}	
equals the value function of our control problem. Recalling Theorem \ref{thm:zerovar}, a straight consequence of equations (\ref{dv2}) and (\ref{cost2go}) therefore is:
\begin{proposition}
The free energy (\ref{cgf}) is equal to
\begin{equation}
\gamma = \bE[Y_0]\,,
\end{equation}
where the expectation is over the initial conditions $X_0$ in $Y_0=V(X_0,0)$.
\end{proposition}

\begin{remark}
A remark on the role of the control variable $Z_s$ in the BSDE is in order. In (\ref{bsde}), let $h=0$ and consider a random variable $\xi$ that is square-integrable and $\mathcal{F}_\tau$-measurable where $\mathcal{F}_s$ is the $\sigma$-Algebra generated by $(B_r)_{0\leqslant r\leqslant s}$. Ignoring the measurability for a second, a pair of processes $(Y,Z)$ satisfying
\begin{equation}\label{bsdeEx}
dY_s = Z_s\cdot dB_s\,,\quad Y_\tau=\xi\,.
\end{equation}
is $(Y,Z)\equiv(\xi,0)$, but then $Y$ is not adapted unless the terminal condition $\xi$ is a.s.~constant, because $Y_t$ for any $t<\tau$ is not measurable \wrt $\mathcal{F}_s\subset \mathcal{F}_\tau$. An adapted version of $Y$ can be obtained by replacing $Y_t=\xi$ by its best approximation in $L^2$, i.e.~by the projection $Y_t=\bE[\xi|\mathcal{F}_t]$. Since the thus defined process $Y$ is a martingale  with respect to our filtration, the martingale representation theorem asserts that $Y_t$ must be of the form
\begin{equation}
Y_t = \bE[\xi] + \int_0^t \tilde{Z}_s\cdot dB_s\,,
\end{equation}
for some unique, predictable process $\tilde{Z}$. Subtracting the last equation from $Y_\tau=\xi$ yields
\begin{equation}
Y_t = \xi - \int_t^\tau \tilde{Z}_s\cdot dB_s\,,
\end{equation}
or, equivalently,
\begin{equation}
dY_t = \tilde{Z}_s\cdot dB_s\,,\quad Y_\tau=\xi\,.
\end{equation}
Hence $Z_s=\tilde{Z}_s$ in (\ref{bsdeEx}) is indeed a control variable that makes $Y$ adapted.
\end{remark}

\begin{remark}
The forward-backward SDE (or: FBSDE) (\ref{fbsde}) is called \emph{uncoupled} since the forward SDE does not depend on the solution to the associated BSDE, a property that will be exploited in various ways later on.
\end{remark}

\subsection{Importance sampling in path space, cont'd.}

The role of the process $Z$ in the FBSDE representation of the dynamic programming equation is not only to guarantee that $Y$ in (\ref{fbsde}) is adapted, so that $Y_t=V(x,t)$ is the value function, but it can be literally interpreted as a control since $Z_t=\sigma(X_t)^T\nabla_xV(X_t,t)$. We could compute the optimal control for the zero-variance importance sampling estimator (\ref{zeroSDE}) by solving (\ref{fbsde}) with initial condition $X_t=X^u_t$ on-the-fly, in which case one has to compute the solution of (\ref{fbsde}) in parallel to the solution of (\ref{drivenSDE2}). Depending on the nature of the system (in particular the state space dimension) this on-the fly-computation, though computationally demanding, may be nonetheless a sensible alternative to numerical schemes that seek to approximate the value function by globally supported basis functions, which may be an ill-conditioned problem, e.g.~if the majority of the trajectories are known to reside inside a small set.

As an alternative that we discuss in detail later on, we suggest to define a feedback control for the controlled SDE (\ref{drivenSDE2}) by
\begin{equation}\label{uCand}
u_t = -\sigma(X^u_t)^T\nabla_x V_K(X^u_t,t)\,,
\end{equation}
where
\begin{equation}\label{VK}
V_K(x,t) = \sum_{k=1}^K \alpha_k(t)\phi_k(x)
\end{equation}
with $\alpha_k\in\R$ and continuously differentiable (e.g.~radial) basis functions $\phi_k$ is an approximation ansatz for the value function. Then, by Girsanov's Theorem,
\begin{equation}\label{IS}
\bE\left[\exp\left(-W_\tau\right)\right] = \bE_Q\!\left[\exp(-L^u_\tau-W^u_\tau)\right]
\end{equation}
where $L^u_\tau=\log(dQ/dP)$ is the log likelihood of the change of measure from $P$ to $Q$ on $\mathcal{F}_\tau$, as given by (\ref{girsanov}). By continuity of the functional (\ref{IS}), we expect that any unbiased estimator of the right hand side of (\ref{IS}) will have a considerably smaller variance than the plain vanilla estimator (based on independent draws from $P$), provided that $V_K\approx V$ approximates the value function.

\section{Least-squares Monte Carlo}\label{sec:mc}

In this section we discuss the numerical discretisation of the uncoupled FBSDE (\ref{fbsde}), following an approach that was first suggested by Gobet et al.~\cite{GobetEtal2005} and later on refined by several authors; here we suggest a semi-parametric approach with radial basis functions based on the work by Bender and Steiner \cite{BenderSteiner2012}.

\subsection{Time stepping scheme}

The fact that the FBSDE (\ref{fbsde}) is decoupled implies that it can be discretised by an explicit time-stepping algorithm. Here we utilise a variant of the least-squares Monte Carlo algorithm proposed in \cite{GobetEtal2005}. The convergence of the numerical schemes for an FBSDE with quadratic nonlinearities in the driver has been analysed in \cite{TurkedjievDissertation2013}. 
The least-squares Monte Carlo scheme is based on the Euler discretisation of (\ref{fbsde}), specifically,
\begin{equation}\label{fbsdeEuler}
\begin{aligned}
\hat{X}_{n+1}  & = \hat{X}_n + \Delta t\, b(\hat{X}_n,t_n) + \sqrt{\Delta t}\,\sigma(\hat{X}_n)\xi_{n+1}\\
\hat{Y}_{n+1} & = \hat{Y}_{n} - \Delta t\, h(\hat{X}_n,\hat{Y}_n,\hat{Z}_n) + \sqrt{\Delta t}\,\hat{Z}_n\cdot\xi_{n+1}\,,
\end{aligned}
\end{equation}
where $(\hat{X}_n,\hat{Y}_n,\hat{Z}_n)$ denotes the numerical discretisation of the joint process $(X_s,Y_s,Z_s)$, where we set $X_s \equiv X_{\tau_O}$ for $s\in(\tau_O,T]$ when $\tau_O<T$, and $(\xi_i)_{i\geqslant 1}$ is an i.i.d.~sequence of normalised Gaussian random variables. Now let
\[
\mathcal{F}_n = \sigma\big(\big\{\hat{B}_k: 0\leqslant k\leqslant n\big\}\big)
\]
be the $\sigma$-algebra generated by the discrete Brownian motion $\hat{B}_n:=\sqrt{\Delta t}\sum_{i\leqslant n}\xi_i$. By definition, the continuous-time process $(X_s,Y_s,Z_s)$ is adapted to the filtration generated by $(B_r)_{0\leqslant r\leqslant s}$. For the discretised process, this implies
\begin{equation}\label{condExp}
\hat{Y}_n = \E\big[\hat{Y}_n|\mathcal{F}_n\big] = \E\big[\hat{Y}_{n+1} + \Delta t \,h(\hat{X}_n,\hat{Y}_n,\hat{Z}_n)|\mathcal{F}_n\big]\,,
\end{equation}
using that $\hat{Z}_n$ is independent of $\xi_{n+1}$.
In order to compute $\hat{Y}_n$ from $\hat{Y}_{n+1}$, it is convenient to replace $(\hat{Y_{n}},\hat{Z}_{n})$ on the right hand side by $(\hat{Y}_{n+1},\hat{Z}_{n+1})$, so that we end up with the fully explicit time stepping scheme
\begin{equation}\label{Yn}
\hat{Y}_n = \E\big[\hat{Y}_{n+1} + \Delta t \,h(\hat{X}_n,\hat{Y}_{n+1},\hat{Z}_{n+1})|\mathcal{F}_n\big]\,.
\end{equation}
Note that we can use the identification of $Z$ with the optimal control (\ref{uCand}) and replace $\hat{Z}_{n+1}$ in the last equation by
\begin{equation}\label{Zn}
\hat{Z}_{n+1} = \sigma(\hat{X}_{n+1})^T \nabla V_K(\hat{X}_{n+1},t_{n+1})\,,
\end{equation}
where $V_K$ is given by the parametric ansatz (\ref{VK}).

\begin{remark}
	If an explicit representation of $\hat{Z}_{n}$ such as (\ref{Zn}) is not available, it is possible to derive a time stepping scheme for $(\hat{Y}_{n},\hat{Z}_{n})$ in the following way: multiplying the second equation in (\ref{fbsdeEuler}) by $\xi_{n+1}\in\R^{m}$ from the left, taking expectations and using the fact that $\hat{Y_{n}}$ is adapted,  it follows that
	\begin{equation}
	0 = \bE\!\left[\xi_{n+1}\big(Y_{n+1} - \sqrt{\Delta t}\hat{Z}_{n}\cdot \xi_{n+1}\big)\big| \mathcal{F}_{n}\right]
	\end{equation}
	or, equivalently,
	\begin{equation}
	\hat{Z}_{n} = \frac{1}{\sqrt{\Delta t}} \bE\!\left[\xi_{n+1}Y_{n+1}\big| \mathcal{F}_{n}\right]\,.
	\end{equation}
	Together with (\ref{Yn}) or, alternatively, with
	\begin{equation}\label{Yn2}
	\hat{Y}_n = \E\big[\hat{Y}_{n+1} + \Delta t \,h(\hat{X}_n,\hat{Y}_{n+1},\hat{Z}_{n})|\mathcal{F}_n\big]\,,
	\end{equation}
	we have a fully explicit scheme for $(\hat{Y}_{n},\hat{Z}_{n})$.
\end{remark}

\subsection{Conditional expectation}

We next address the question how to compute the conditional expectations with respect to $\mathcal{F}_{n}$. To this end, we recall that the conditional expectation can be characterised as a best approximation in $L^{2}$:
\[
\E\big[S|\mathcal{F}_n\big] = \mathop{\rm argmin}_{Y\in L^2,\, \mathcal{F}_n \textrm{-measurable}}\E[|Y-S|^2]\,.
\]
(Hence the name \emph{least-squares Monte Carlo}.)
Here measurability \wrt $\mathcal{F}_{n}$ means that $(\hat{Y}_{n},\hat{Z}_{n})$ can be expressed as functions of $\hat{X}_{n}$. In view of the ansatz (\ref{VK}) and equation (\ref{Yn}), this suggests the  approximation scheme
\begin{equation}\label{condVar}
\hat{Y}_n \approx \mathop{\rm argmin}_{Y=Y(\hat{X}_n)}\frac{1}{M}\sum_{m=1}^{M}\left|Y -  \hat{Y}_{n+1}^{(m)} - \Delta t \, h\big(\hat{X}^{(m)}_n,\hat{Y}^{(m)}_{n+1},\hat{Z}^{(m)}_{n+1}\big)\right|^2\,,
\end{equation}
where the data at time $t_{n+1}$  is given in form of $M$ independent realisations of the forward process, $\hat{X}_n^{(m)}$, $m=1,\ldots,M$, the resulting values for $\hat{Y}_{n+1}$,
\begin{equation}\label{Yn+1}
\hat{Y}^{(m)}_{n+1} = \sum_{k=1}^K \alpha_k(t_{n+1})\phi_k\big(\hat{X}^{(m)}_{n+1}\big)\,,
\end{equation}
and
\begin{equation}\label{Zn+1}
\hat{Z}^{(m)}_{n+1} = \sigma\big(\hat{X}^{(m)}_{n+1}\big)^T\sum_{k=1}^K \alpha_k(t_{n+1})\nabla\phi_k\big(\hat{X}^{(m)}_{n+1}\big)\,.
\end{equation}
At time $T:=N\Delta t$, the data are determined by the terminal cost:
\begin{equation}\label{YZterm}
\hat{Y}^{(m)}_N = g\big(X^{(m)}_N\big)\,,\quad \hat{Z}^{(m)}_N = \sigma\big(\hat{X}^{(m)}_{N}\big)^T\nabla g\big(X^{(m)}_N\big)
\end{equation}
Note that we have defined the forward process so that all trajectories have length $T$, but the realisations may be constant between $\tau_O$ and the terminal time $T$.

The unknowns that have to be computed in every iteration step are the coefficients $\alpha_k$, which makes them functions of time, i.e.~$\alpha_k=\alpha_k(t_{n+1})$.
We call $\hat{\alpha}=(\alpha_1,\ldots,\alpha_K)$ the vector of the unknowns, so that the least-squares problem that has to be solved in the $n$-th step of the backward iteration is of the form
\begin{equation}\label{leastSq}
\hat{\alpha}(t_n) = \mathop{\rm argmin}_{\alpha\in\R^K} \left\|A_n\alpha - b_n\right\|^2\,,
\end{equation}
with coefficients
\begin{equation}\label{leastSqA}
A_n = \left(\phi_k\big(\hat{X}_n^{(m)}\big)\right)_{m=1,\ldots,M;k=1,\ldots,K}\,
\end{equation}
and data
\begin{equation}\label{leastSqb}
b_n = \left(\hat{Y}_{n+1}^{(m)} + \Delta t\,h\big(\hat{X}^{(m)}_n,\hat{Y}^{(m)}_{n+1},\hat{Z}^{(m)}_{n+1}\big)\right)_{m=1,\ldots,M}\,.
\end{equation}
Assuming that the coefficient matrix $A_n\in\R^{M\times K}$, $K\leqslant M$ defined by (\ref{leastSqA}) has maximum rank $K$, then the solution to (\ref{leastSq}) is given by
\begin{equation}\label{leastSqSol}
\hat{\alpha}(t_n) = \left(A_n^T A^{}_n\right)^{-1}A_n^T b^{}_n\,.
\end{equation}

\begin{algorithm}
	\caption{Least-squares Monte Carlo}\label{lsmc}
	\begin{spacing}{1.1}
		\begin{algorithmic}
			\State Define $K,M,N$ and $\Delta t = T/M$.
			\State Set initial condition $x\in\R^d$.
			\State Choose radial basis functions $\{\phi_k\in C^1(\R^d,\R)\colon k=1,\ldots,K\}$.
			\State Generate $M$ independent realisations $\hat{X}^{(1)},\ldots,\hat{X}^{(M)}$ of length $N$ from
			\[
			\hat{X}_{n+1} = \hat{X}_n + \Delta t\, b(\hat{X}_n,t_n) + \sqrt{\Delta t}\,\sigma(\hat{X}_n)\xi_{n+1}\,,\; \hat{X}_0 = x\,.
			\]
			\State Initialise BSDE by
			\[
			\hat{Y}^{(m)}_{N} =  g\big(\hat{X}^{(m)}_{N}\big)\,,\quad \hat{Z}^{(m)}_{N} =  \sigma\big(\hat{X}^{(m)}_{N}\big)^T\nabla g\big(\hat{X}^{(m)}_{N}\big)\,.
			\]
			\For {$n=N-1\colon 1$}
			\State Assemble linear system $A_n\hat{\alpha}(t_n)=b_n$ according to (\ref{leastSq})--(\ref{leastSqb}).
			\State Evaluate $\hat{Y}^{(m)}_{n}$ and $\hat{Z}^{(m)}_{n}$ according to
			\[
			\hat{Y}^{(m)}_{n} =  \sum_{k=1}^K \alpha_k(t_{n})\phi_k\big(\hat{X}^{(m)}_{n}\big)\,,\quad \hat{Z}^{(m)}_{n} =  \sigma\big(\hat{X}^{(m)}_{n}\big)^T\sum_{k=1}^K \alpha_k(t_{n})\nabla\phi_k\big(\hat{X}^{(k)}_{n}\big)\,.
			\]
			\State If necessary, adapt basis functions $\phi_k$.
			\EndFor
		\end{algorithmic}
	\end{spacing}
\end{algorithm}

The thus defined scheme that is summarised in Algorithm \ref{lsmc} is strongly convergent of order 1/2 as $\Delta t\to 0$ and $M,K\to\infty$; see \cite{GobetEtal2005}. Controlling the approximation quality for finite values $\Delta t, M, K$, however, requires a careful adjustment of the simulation parameters and appropriate basis functions, especially with regard to the condition number of the matrix $A_n$, and we will discuss suitable strategies to determine a good basis in the next section.

\begin{remark}\label{Modif}
The accuracy of the solution to the backward SDE depends on wether the distribution of the terminal condition $g(X_{\tau})$ is accurately sampled. If the forward process is metastable, however, it may happen that  $g(X_{\tau})$ is poorly sampled. In this case, it is possible to change the drift of the forward SDE from $b$ to, say, $b_0$ where $b_0$ is chosen such that the forward trajectories densely sample the statistic $g(X_{\tau})$, without affecting the value function or the resulting optimal control: Assuming that the noise coefficient $\sigma$ is square and invertible, it is easy to see that the dynamic programming PDE (\ref{hjb}) can be recast as
\begin{align*}
\partial_{t} V + \tilde{L}V + \tilde{h}(s,x,V,\sigma^{T}\nabla_{x}V) & = 0\quad \textrm{in $E$}\\ V & = g\quad \textrm{on $\partial E^{+}$}\,,
\end{align*}
where
\[
\tilde{L} = L  - (b-b_0)\cdot\nabla
\]
is the generator of a forward SDE with drift $b_0$, and
\[
\tilde{h}(x,y,z) = h(x,y,z) +  \sigma(x)^{-1}(b(x)-b_0(x))\cdot z\,
\]
is the driver of the corresponding backward SDE. Hence we can change the drift of the forward SDE at the expense of modifying the running cost, without affecting the optimal control. Changing the drift may be moreover advantageous in connection with the martingale basis approach of Bender and Steiner \cite{BenderSteiner2012} who have suggested to use basis functions that are defined as conditional expectations of certain linearly independent candidate functions over the forward process, which makes the basis functions martingales. Computing the martingale basis, however, comes with a large computational overhead, which is why the authors consider only cases in which the conditional expectations can be computed analytically. Changing the drift of the forward SDE may thus be used to simplify the forward dynamics so that its distribution becomes analytically tractable.	
\end{remark}

\section{Numerical illustration}\label{sec:num}

We shall illustrate the previous considerations with a standard example. To this end, we consider a one dimensional diffusion in the double-well potential $U(x)=(x^2-1)^2$ that is governed by the equation
 \begin{equation}\label{sdegrad1}
dX_{s} = -\nabla U(X_{s})ds + \sigma dB_{s} \,,\quad X_0=x\,,
\end{equation}
and want to compute the probability of exiting from the left well $O = \left\{ x < 0\right\}$ before time $T < \infty$.
More specifically, we  set $f\equiv 0$ and $g(x)= - \log\left(\mathbf{1}_{\partial O}(x)\right)$ in equation (\ref{Wfunc}) and define the bounded stopping time $\tau=\tau_{O}\wedge T$ to be the minimum of the first exit time $\tau_{O}$ of the set $O$ and the terminal time $T$. Note that $\tau_{O}$ is a.s.~finite since the potential $U$ is growing sufficiently fast at infinity, so that $(X_{s})_{s\geqslant 0}$ is Harris recurrent.  

For the equivalent stochastic control problem with the cost
\begin{equation}
J(u) = \bE\!\left[\frac{1}{2}\int_{0}^{\tau} |u_s|^2\,ds - \log\left(\mathbf{1}_{\partial O}(X^u_\tau)\right)\right]
\end{equation}
and the controlled process
 \begin{equation}
dX^{u}_{s} = \left(\sigma u_{s} - \nabla U(X^{u}_{s})\right)ds + \sigma dB_{s} \,,\quad X^{u}_0=x\,,
\end{equation}
this means that the control $u$ seeks to push the process towards the set boundary $\partial O$ when $s\approx T$ and the process has not yet left the set $O$, for otherwise there will an infinite cost to pay. 

Since such an infinite terminal cost is numerically difficult to handle, we consider a regularised control problem and replace $g$ by $g^{\varepsilon}= -\log(\mathbf{1}_{\partial O}(x) + \varepsilon)$; for the numerical calculations, we choose $\varepsilon = 0.01$. The duality relation (\ref{dv}) between the control value  $\gamma^{\varepsilon}=\min_{u} J(u)$ for fixed  initial data $X_{0}=x$ and the transition probability $P(\tau_{O}<T)$ then reads 
\begin{equation}
P(\tau_{O}<T | X_{0} = x) = \exp(-\gamma^{\varepsilon}) - \varepsilon\,. 
\end{equation}

We will compare the results from the FBSDE solution for $\gamma^{\varepsilon}$ with a reference solution that is obtained from numerically solving the linear PDE
\begin{equation}
 \left(\frac{\partial}{\partial t} - L\right) \psi(x,t) = 0 \,,\quad (x,t)\in O\times [0,T)
\end{equation}
together with the boundary conditions\footnote{For the numerical computation, we add reflecting boundary conditions at $x=-L$ for some $L>0$, the precise value of which does not affect the results (assuming that it is sufficiently large, say, $L>3$) since the potential has a 4-th order growth.}
\begin{equation}
\begin{aligned}
\psi(0,t) & = 1\,, \quad t\in [0,T)\\ 
\psi(x,0) & = 0\,, \quad x\in O\,.
\end{aligned}
\end{equation}
Then 
\begin{equation}
\psi(x,T) = P(\tau_{O}<T | X_{0} = x)\,.
\end{equation}

Table \ref{numres} below shows the reference value $V^{\varepsilon}_{ref}(0,x) := - \log \left( \psi(x,T) + \varepsilon\right)$, together with the corresponding FBSDE solution.
The procedure to obtain the FBSDE solution is described in Algorithm \ref{lsmc}, and the table displays the results for different values of $K,M,N=\lfloor T/\Delta t\rfloor$. As basis functions we choose 
\begin{equation}
\phi^{\mu_k, \delta}_{k,n}(x) = \exp\left(-\frac{(\mu_k-x)^2}{2 \delta}\right)\,,
\end{equation}
where $\delta = 0.1$ is fixed but $\mu_k = \mu_k(n)$ varies with time such that the forward process can be well covered by the basis functions. More precisely, the centres of the basis functions are chosen by simulating $K$ additional independent forward trajectories $X^{(k)}, k=1,\ldots,K $ and letting $\mu_k(n) = X^{(k)}_n$. We let the whole algorithm run 20 times and compute empirical mean and variance of $V^{\varepsilon}$, denoted by $\bar{V}^{\varepsilon}$ and $ S^2(V^{\varepsilon})$. The results are shown in the table.

\begin{center}
\captionof{table}{Numerical results for the FBSDE scheme described in Algorithm 1.} \label{numres}
    \begin{tabular}{ p{5.5cm}  || c | c | c }
 & $\ V_{ref}^{\varepsilon}(0,x)\ $& $\ \bar{V}^{\varepsilon}(0,x)\ $  & $\ S^2(V^{\varepsilon}(0,x))\ $  \\     \hline \hline
 $K=8,$ $M = 300,$ $T=5,$ $\Delta t = 10^{-3},$ $x=-1$, $\sigma =1 $ & 0.3949 &0.3748 &$10^{-3}$  \\ \hline
$K=5,$ $M = 300,$ $T=1,$ $\Delta t = 10^{-3},$ $x=-1$, $\sigma =1 $ & 1.7450 &1.6446 &0.0248  \\ \hline
 $K=5,$ $M = 400,$ $T=1,$ $\Delta t = 10^{-4},$ $x=-1,$ $\sigma=0.6$ & 4.3030 & 4.5779 & $10^{-3}$ \\ \hline
 $K=6,$ $M = 450,$ $T=1,$ $\Delta t = 10^{-4},$ $x=-1,$ $\sigma=0.5$ & 4.5793 & 4.6044 & $5\cdot 10^{-4} $
   \end{tabular}
\end{center}

Overall we find that the FBSDE scheme results in a fairly good approximation of the value function and, as a consequence of the smoothness of the basis functions, of the optimal control. Moreover, due to the adaptive choice of the basis functions $\{\phi^{\mu_{k},\delta}\}$, the results do not seem to be very sensitive to the noise intensity $\sigma$ or the time horizon $T$. Speaking of which, we stress that increasing the number of basis functions $K$ is not always advisable, since the  matrix $A$ in \eqref{leastSqA} can easily become rank deficient, especially if $\sigma $ is small and the trajectories stay close together. Therefore it is crucial to check the rank of $A$ in the simulation and to set $K$ to the value of the maximally observed rank. 

\subsection{Computational issues}

Let us also discuss the fact that we set $X_s \equiv X_{\tau_O}$ for $s\in(\tau_O,T]$ when $\tau_O<T$ again in more detail. Setting the forward trajectories constant from the exit time on, allows to include the terminal condition $g(X_{\tau \wedge T})$ into the least squares problem at time $T$, i.e. into the initialisation step $b_N$, for all backward trajectories. It seems that this stabilises the solution of the backward trajectory $\hat{Y}$. Another approach, following the equations more closely, would be to start each backward trajectories individually from either $\tau$ or $T$ depending on whether the corresponding forward trajectory $\hat{X}$ has made an exit or not. This approach induces numerical problems, though, because the data vector (\ref{leastSqb})---that would normally be dominated by the positive term  $\hat{Y}_{n+1}$ when all backward trajectories were starting from $T$---is now perturbed at the different exit times by the negative value $-\log(\varepsilon)$. This renders the solution $\hat{\alpha}_n$ of the linear equation \eqref{leastSqSol} rougher, which in turn leads to fluctuations in the solution of $\hat{Y}_n$ and $\hat{Z}_n$ which can build up and eventually lead to an explosion of the solutions. 

Let us further make suggestions how to efficiently treat the case when $T$ is large. We will resort to the ideas of Remark \ref{Modif} here, which suggests to modify the drift $b$ to $b_0$ such that under the new drift the event which determines the stopping time $\tau$ is not rare anymore. Assume now, that for all trajectories $\hat{X}^{(m)}, m=1, \ldots, M$ the family of stopping times
\begin{equation}
\tau_O^m = \left\{s > 0: X^{(m)}_s \notin O \right\}
\end{equation}
is dominated by $T$ in the sense that 
\begin{equation}
\tilde{T} := \max\{\tau^m_O\colon m=1,\ldots, M\}  \ll T\,.
\end{equation}
Then the terminal condition $g$ is essentially known at time $\widetilde{T}$ and the same is true for the backward dynamics. Hence, we suggest in case that $T$ is large to modify the drift such that $\widetilde{T}$ will be small and run the algorithm only up to time $\tilde{T}.$ In this case we propose to start each backward trajectory individually from the corresponding exit time on. The matrix $A_n$ is then of size $K \times M_n$ where 
\begin{equation}
M_n = \left|\left\{m:\hat{X}_{n-1}^{(m)} \in O\right\}\right| 
\end{equation}
is the number of trajectories which have not left the set $O$ up to time step $n$. This ensures that $A$ is not rank deficient at these times which would be the case if we set all trajectory constant after the exit, due to the definition of $A_n$ with 
\begin{equation}
\left(A_n\right)_{k,m} = \phi^{\mu_k, \delta}_{k,n}(\hat{X}^{(m)}_n)
\end{equation}
 because the basis functions are evaluated at the same constant value for all these trajectories. 
To the best of our knowledge, the approximation error of the least squares Monte Carlo algorithms with random stopping  times has not been analysed so far, and we leave this topic for future work.  

We want to add that in contrast to the complexity of numerically solving the HJB equation, which grows exponentially in the dimension $d$, the complexity of solving the FBSDE is determined by solving the SDE and linear equations, i.e. is at most cubic in $d$ and in the number $K$ of basis functions.

\section{Conclusion and outlook}\label{sec:out}
We have presented a numerical method to compute the free energy of path space functionals of a diffusion process where the functionals may depend on paths having a random length. Free energies of path space functionals appear in connection with rare event simulation and, as a guiding example for this article, we have considered exit probabilities that are relevant in the context of molecular dynamics or risk analysis. 

The approach for efficiently computing path space free energies is based on a variational characterisation of the free energy as the value function of an optimal control problem or, equivalently, as an adaptive importance sampling strategy that is based on the optimal control of the aforementioned stochastic control problem; 
as we have argued, the importance sampling estimator for the free energy enjoys a minimum variance property under the optimal control. 
Our numerical strategy for solving the underlying stochastic control problem is based on the reformulation of the corresponding semi-linear dynamic programming equation as a forward-backward stochastic differential equation, which can be solved quite efficiently using a least squares Monte Carlo method. For our guiding example, the reformulation of the adaptive importance sampling algorithm as a forward-backward SDE showed promising results.

We have discussed several options that can help to improve the convergence of the least squares algorithm. For example, we have discussed the option of changing the drift of the forward SDE by modifying the cost functional of the corresponding control problem; while this does not change the dynamic programming equation of the underlying control problem, the corresponding forward-backward stochastic differential equations are different, and it is possible to control the speed of convergence of the numerical method in this way, by controlling the random length of the forward trajectories.  

Another aspect that we have only briefly touched upon is the choice of the basis functions for the least squares algorithm. A convenient choice are martingale basis functions that, by definition, are non-parametric and adaptive. Evaluating the martingale basis requires to compute on-the-fly conditional expectations and it is possible to change the drift of the forward SDE so as to avoid numerically expensive computations of the conditional expectations. In this article we used a semi-parametric approach, and future research should address the non-parametric one.
Another interesting topic concerns sampling problems on an infinite time horizon, which can be represented by a stopping time for hitting an \emph{impossible set}, a set which the dynamics can never reach. 

We believe that forward-backward SDE are an interesting numerical and analytical tool for applications in computational statistical mechanics that connects such diverse topics as control, filtering and estimation. A specific feature of the proposed method is that the corresponding forward-backward SDE are decoupled, which leaves room for combining the aforementioned tasks with coarse-graining and model reduction techniques. We leave all this for future work.

\subsubsection{Acknowledgement}
This research has been partially funded by Deutsche For\-schungsgemeinschaft (DFG) through the grant CRC 1114 "Scaling Cascades
in Complex Systems", Project A05 "Probing scales in equilibrated systems by optimal nonequilibrium forcing". Omar Kebiri received funding from the EU-METALIC II Programme.

\end{document}